\documentclass[11pt,twoside]{article}
\usepackage{latexsym,amssymb,amsfonts}
\usepackage{amsmath}
\usepackage{amscd}
\usepackage[english]{babel}
\usepackage{pstricks}
\usepackage{pst-node}
\usepackage{multido}
\def\co{\colon\thinspace}

\begin{document}

\title{Kauffman state sums and bracket deformation}
\author{Nikos Apostolakis\footnote{partially supported by PSC CUNY 
    Research Award PSCOOC-37-105 }\\ 
  \begin{tabular}{c}
    BCC of CUNY\\
\texttt{nikolaos.apostolakis@bcc.cuny.edu}
  \end{tabular}
 \and Uwe Kaiser \footnote{partially
supported by NSF grant DMS-0204627}\\ 
  \begin{tabular}{c}
    Boise State University\\
  \texttt{kaiser@math.boisestate.edu}
\end{tabular}}

\date{}

\maketitle

\abstract{\noindent We derive a formula expanding the bracket with
respect to a natural deformation parameter. The expansion is in terms
of a two-variable polynomial algebra of diagram resolutions generated by
basic operations involving the Goldman bracket. A functorial
characterization of this algebra is given. Differentiability
properties of the star product underlying the Kauffman bracket are
discussed.}

\vskip .2in


\vskip .1in

\noindent \textsf{Mathsubject Class.:} 57M25, 57M35, 57R42

\vskip .5in

\section{Introduction.}

Throughout let $F$ be a compact connected oriented surface, $\amalg$
means disjoint union. 

\vskip .1in

Let $\mathfrak{C}(F)$ be the set of isotopy classes of closed
$1$-dimensional submanifolds of the interior of $F$ without
inessential components. These are called \textit{curve systems} on
$F$. There is the empty curve system $\emptyset\in \mathfrak{C}(F)$.

Let $\mathfrak{L}(F)$ denote the set of isotopy classes of framed
unoriented links in $F\times I$, including the empty link
$\emptyset$. The set $\mathfrak{L}(F)$ is identified
with the set of isotopy classes of diagrams $\mathfrak{D}(F)$
on $F$ up to Reidemeister moves of type $II$ and $III$ (diagram will
always mean regular diagram).
The identification is given by regular projection and blackboard framing.
 
Let $\mathsf{k}$ be a commutative ring with $1$.
It is a result of Przytycki \cite{P} that $\mathfrak{C}(F)$ is a module
basis of the Kauffman bracket skein module $\mathfrak{K}(F)$ of $F$.
By definition, $\mathfrak{K}(F)$ is the quotient of the free module 
$\mathsf{k}[t,t^{-1}]\mathfrak{L}(F)$ by the submodule generated by
the elements 
$K_++tK_{0}+t^{-1}K_{\infty}$ (resolution) 
and $K\amalg U+(t^2+t^{-2})K$ (trivial component).

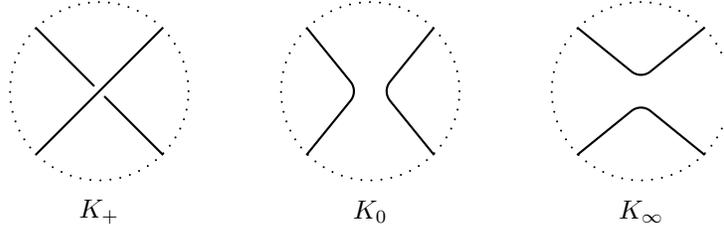
\begin{center}
\unitlength 0.8cm
\psset{unit=\unitlength}
  \begin{pspicture}(-2,-2)(11,2)
    \pscircle[linestyle=dotted](0,0){1.5}
    \psline(1.061,-1.061)(-1.061,1.061)
    \psline[border=.1](1.061,1.061)(-1.061,-1.061)
    \rput(0,-2){$K_{+}$} \rput(4.5,0){%
      \pscircle[linestyle=dotted](0,0){1.5}
      \psline[linearc=.25](-1.061,-1.061)(-.2,0)(-1.061,1.061)
      \psline[linearc=.25](1.061,-1.061)(.2,0)(1.061,1.061) }%
    \rput(4.5,-2){$K_{0}$} \rput(9,0){%
      \pscircle[linestyle=dotted](0,0){1.5}
      \psline[linearc=.25](-1.061,1.061)(0,.2)(1.061,1.061)
      \psline[linearc=.25](-1.061,-1.061)(0,-.2)(1.061,-1.061) }%
    \rput(9,-2){$K_{\infty}$}
  \end{pspicture}
\end{center}
 
The over-under crossing information of $K_+$
and the orientation of $F$ determine the resolutions $K_0$ and
$K_{\infty}$ in the usual way. $U$ denotes a component whose
projection is an embedded circle on $F$, which bounds a disk in the
complement of the projection of $K_+$. 

\vskip .1in
 
Let $[K]$ respectively $[D]$ denote the image of a framed link
respectively diagram in the 
Kauffman bracket module.     
The isomorphism 
$$\mathfrak{K}(F)\rightarrow \mathsf{k}[t,t^{-1}]\mathfrak{C}(F)$$
is established using the Kauffman bracket state sum of diagrams
$$<D>=\sum_{\sigma}(-t)^{\zeta(\sigma)-\iota(\sigma )}(-t^2-t^{-2})^{\mu( \sigma)}D(\sigma )\in \mathsf{k}[t,t^{-1}]\mathfrak{C}(F)$$  
for each diagram $D$ on $F$.
The sum is over all Kauffman states $\sigma $ of $D$, i.e.
assignments of state markers $0,\infty$ to each crossing of the diagram.
The functions $\zeta$, respectively $\iota$, assign to each state  
its number of $0$-, respectively $\infty$-markers.
Recall that the assignment of a state marker to a crossing defines a resolution
of that crossing. Then $\mu $ assigns to the state $\sigma $ the number 
of inessential circles in the resolution determined by $\sigma $. 
$D(\sigma )$ is the collection of essential components which 
appear in the resolution determined by $\sigma $.
Note that  $(\zeta +\iota )(\sigma )$ is equal to the 
number $c$ of crossings of $D$ for all $\sigma $.

\vskip .1in
 
It is easy to see that $<D>$ only depends on $[D]$. 
Using the inclusion $\mathfrak{C }(F)\subset \mathfrak{D}(F)$ it
follows that
$[D]\mapsto <D>$ defines a module isomorphism.   

\vskip .1in

Recall that the module $\mathfrak{K}(F)$ actually is a
$\mathsf{k}[t,t^{-1}]$-algebra with multiplication $\star$ defined by 
stacking links.
For two framed links $K,K'$ in $F\times I$ we let $[K]\star [K']$ be the 
element of $\mathfrak{K}(F)$ represented by placing $K\subset F\times
[1,2]$ and $K'\subset F\times [0,1]$,
thus $K\amalg K'\subset F\times [0,2]$ (and 
$[0,2]\cong [0,1]$ in a natural way). 
For two diagrams $D,D'$ on $F$ we let $D\triangleright D'$ denote \textit{some} diagram on $F$ defined by having only crossings of $D$ over $D'$. Note that 
$D\triangleright D'$ is \textit{not} a well-defined diagram but $<D\triangleright D'>\in \mathsf{k}[t,t^{-1}]\mathfrak{C}(F)$ is well-defined because of the isotopy invariance of the Kauffman bracket. For 
$\alpha ,\beta \in \mathfrak{C}(F)$ let $\alpha \star \beta $ denote the 
result of multiplication in $\mathfrak{K}(F)$, and expanding using the Kauffman bracket.
Thus $\alpha  \star \beta \in \mathsf{k}[t,t^{-1}]\mathfrak{C}(F)$. 
Note that the $\star $-product is non-commutative except for $F$ a disk, annulus or $2$-sphere. 

It is a difficult problem to relate the expansions $<D>$ respectively
$\alpha \star \beta $  with the geometry of the diagram respectively
curves. For general diagrams this is a nontrivial question even in the
commutative case. The $\star $-product of two curve systems
on $F$ is trivially known in the commutative case but its computation
is difficult in the non-commutative case. A complete answer indicating 
the relation of this problem with non-commutative geometry has been
given for $F=S^1\times S^1$ by Frohman and Gelca \cite{FG}.  
 
It is our goal to study the 
combinatorics of the \textit{deformation theory} of the Kauffman
bracket and the $\star $-product. Assume for the moment that
$\mathsf{k}$ is a field of characteristic $0$ (not necessarily
algebraically closed). Then there is an embedding 
$\mathsf{k}[t,t^{-1}]\rightarrow \mathsf{k}[[h]]$ defined by mapping
$t$ to $e^h$. 
Using the inclusions 
$$\mathsf{k}[t,t^{-1}]\mathfrak{C}(F)\rightarrow \mathsf{k}[[h]]\mathfrak{C}(F)\subset 
\mathsf{k}\mathfrak{C}(F)[[h]]$$
we can map $<D>$ into $\mathsf{k}\mathfrak{C}(F)[[h]]$.
(The second inclusion is proper because $\mathfrak{C}(F)$ is an infinite set.) 
The image of the bracket in $\mathsf{k}\mathfrak{C}(F)[[h]]$ is still
denoted $<\ >$ and we can write
$$<D>=\sum_{j=0}^{\infty}<D>_jh^j$$
with $<D>_j\in \mathsf{k}\mathfrak{C}(F)$.   
In the case of $\alpha ,\beta \in \mathfrak{C}(F)$ this defines
$$\alpha \star \beta =\sum_{j=0}^{\infty}\lambda_j(\alpha ,\beta)h^j$$
where for $j\geq 0$ the $\lambda_j$ extend to $\mathsf{k}$-bilinear
mappings
$$\lambda_j\co \mathsf{k}\mathfrak{C}(F)\otimes
\mathsf{k}\mathfrak{C}(F)\rightarrow \mathsf{k}\mathfrak{C}(F).$$ 
Note that this sequence determines the $\star $-product.

For a given diagram $D$ the contribution $<D>_0$ of the state sum can
 be calculated by applying the skein relations for $t=1$ thus
$K_++K_0+K_{\infty}=0$ and $U+2=0$ to the diagram. 
It has been shown by Bullock, Frohman and Kania-Bartoszynska
\cite{BFK} that
$$<D>_1=\sum_{\textrm{crossings}\ p}<D_{p,\infty}>_0-<D_{p,0}>_0,$$
where $D_{p,0}$ respectively $D_{p,\infty}$ are the diagrams resulting 
from the $0$- respectively $\infty$-resolution of the crossing $p$. 
In fact in \cite{BFK} the formula is only given for 
the first order contribution $\lambda_1(\alpha,\beta )$ of the 
$\star $-product of two simple closed curves. But it is easy to see 
that their combinatorial argument 
immediately applies to all diagrams.

The interest in the two results above comes from its relation with
the representation theory of the fundamental group $\pi_1(F)$ of the
surface $F$. In fact let $Rep(F)$ denote the universal
$SL(2,\mathsf{k})$-character ring of $\pi_1(F)$.  It has been shown
by Bullock \cite{B} and Przytycki-Sikora \cite{PS} that $Rep(F)$
is naturally isomorphic with the $\mathsf{k}$-algebra structure on
$\mathsf{k}\mathfrak{C}(F)$ defined by $\lambda_0$. In the case of
an algebraically closed field $\mathsf{k}=\mathsf{K}$ the algebra
$\mathsf{K}\mathfrak{C}(F)$ can be identified with the ring of
character functions $\mathcal{A}(F)$, i.\ e.\ regular functions on
the variety of $SL(2,\mathsf{K})$-representations that are defined
by evaluations and taking traces.  More precisely the isomorphism of
$\mathsf{K}$-algebras is given by
$$\mathfrak{C}(F)\ni \alpha \mapsto n_{\alpha}\in \mathcal{A}(F)$$
where $n_{\alpha }(\rho ):=-tr(\rho (\alpha ))$ for each
representation 
$\rho \co\pi \rightarrow SL(2,\mathsf{K})$. 

It is the main result of \cite{BFK} that for $\mathsf{K}=\mathbb{C}$, 
$$\lambda_1(\alpha,\beta )=\{n_{\alpha} ,n_{\beta } \},$$
where $\{\ ,\ \}$ is the Poisson bracket on $\mathcal{A}(F)$ defined 
from the complex symmetric bilinear form 
$$B(x,y)=-\frac{1}{2}tr(xy)$$
on the Lie algebra $sl(2,\mathbb{C})$ following Goldman \cite{G}.
Recall that for \textit{closed} surfaces this Poisson bracket is defined
from a complex symplectic structure derived from Poincare duality on
$F$ and 
$B$ (see \cite{G2}).

In fact Bullock, Frohman and Kania-Bartoszynski prove that 
the $\star $-product on $\mathbb{C}\mathfrak{C}(F)[[h]]$ as above 
defines a deformation of the algebra $\mathcal{A}(F)$ in the sense of deformation quantization. (In \cite{BFK} it is also shown that
$\mathbb{C}\mathfrak{C}(F)[[h]]$ is isomorphic to a completed Kauffman bracket
algebra $\hat{\mathfrak{K}}(F)$ defined from $\mathbb{C}\mathfrak{L}[[h]]$ 
dividing out by the closure of the submodule defined from the skein relations as before using the substitutions $t=e^h$. We only like to point out 
that all results extend to the case of an algebraically closed field 
$\mathsf{K}$.   

It is our goal in this paper to prove the following result generalizing 
the combinatorial first order formula of \cite{BFK}.

\vskip .1in

\noindent \textbf{Theorem (nontechnical version).} \textit{For all
$j\geq 0$ and rings $\mathsf{k}$ of characteristic $0$ the $j$-th order term 
$$<D>_j\in \mathsf{k}\mathfrak{C}(F)$$
is a sum of diagram resolutions of order $\leq j$, which are 
combinatorial generalizations of $<\ >_0$ and $<\ >_1$.
Corresponding statements hold for the pairings $\lambda_j$.}

\vskip .1in

Note that for all $j\geq 0$, $<\ >_j$ is invariant under Reidemeister
moves of type II and III. In section 5 we will actually prove a
relative version of theorem 1 for diagrams which possibly contain
proper arcs.

\vskip .1in

The result may seem surprising at first. But if one thinks about the
identification of curve systems with regular functions, and observes
that resolutions formally behave like derivatives, it could be
expected, because regular functions are restrictions of polynomials.

\vskip .1in

The non-technical statement above will have to be refined in the
following. In particular we will define the notion of
\textit{order}.  Roughly a diagram resolution of order $j$ is
defined by state summations over $j$-element subsets of the the set
of crossings with state-contributions depending only on the number
of $\infty$-states of the state, followed by state-summation over
the remaining crossings with contributions determined by the number
of trivial components weighted with coefficients in $\mathsf{k}$
(obtained by expanding the unknot contribution in powers of $h$).

\vskip .1in

In section 6 we will apply the theorem above to study
differentiability properties of the deformation. It turns out that the
loop correction terms in the bracket imply that the
bilinear maps $\lambda_j$ are \textit{not} differential
operators of order $\leq j$ in the usual sense. We will discuss 
a combinatorial version of differentiability in section 6.

\vskip .1in

\noindent \textbf{Problem.} Interpret the higher order terms $\lambda_j$ 
for $j\geq 2$ in terms of the geometry of the character variety.

\vskip .1in

\noindent \textbf{Remark.} In \cite{AMR} it is shown that the Poisson
bracket on $\mathcal{A}(F)$ is inherited from a Poisson algebra of
chord diagrams. Moreover for surfaces with nonempty boundary the 
bracket deformation is inherited from the Kontsevitch integral as
constructed in \cite{AMR}. Interestingly in this case the
representation variety has no global symplectic structure. Thus 
it seems particularly interesting to get a better understanding of the
situation for closed surfaces.   
 
\vskip .1in
 
In section 2 we will discuss invariance properties of the bracket
 deformation, in particular invariance under the mapping class group. 
This is to motivate our approach to the diagram resolution algebra
in section 3, where we first give an \textit{intrinsic} description
from functorial properties.  
In section 4 we will prove existence and identify this algebra
with a polynomial algebra in two variables.
In section 5 we prove the technical version of the theorem.
In section 6 we discuss the question of differentiability of the 
star product defined by the Kauffman bracket.

\vskip .1in

The authors would like to thank the referee for pointing out
several typographical mistakes and weak formulations.

\vskip .5in

\section{Invariance properties of the bracket deformation.}

The mapping class group $\mathcal{M}(F)=\pi_0(Diff^+(F))$ of the
surface $F$ acts on the set of curve systems $\mathfrak{C}(F)$ in a
natural way. This action obviously extends to the Kauffman bracket
algebra and is compatible with the multiplication, that is 

$$g(\alpha \star \beta )=(g\alpha )\star (g\beta )$$
for $\alpha, \beta \in \mathfrak{C}(F)$ and $g\in \mathcal{M}(F)$.  
Note that in the calculation of bracket of a diagram $D$ on $F$ we have
$$<gD>=g<D>=\sum_{\sigma }(-t)^{(\zeta -\iota )(\sigma)}(-t^2-t^{-2})^{\mu(\sigma )}(g\delta(\sigma ))$$
reducing the action on diagrams to the action on curve systems. 
Note that $\mathcal{M}$ acts trivially on inessential components.  
The observation above means that the deformation of the commutative product on $\mathcal{A}(F)$ defined by the bracket is invariant under the 
action of the mapping class group. It has been observed by Goldman that
the symplectic structure and thus the Poisson bracket are invariant under the action of $\mathcal{M}(F)$. 

The ring $\mathsf{k}[t,t^{-1}]$ has the natural involution defined
by $t\mapsto t^{-1}$. This defines an anti-involution of the module
$\mathsf{k}[t,t^{-1}]\mathfrak{C}(F)$. Let $\tau
\co\mathfrak{K}(F)\to \mathfrak{K}(F)$ be the anti-involution
defined by changing all crossings of diagrams.  This can be
interpreted as the action of the element of $\pi_0(Diff(F\times
I))$, which is defined by the reflection $I\times I, t\mapsto
(1-t)$. It follows immediately from the definition of the bracket
that $\tau $ corresponds to the ring involution under the bracket
isomorphism, i.e.
$$<\tau (D)>=\tau <D>$$
for all diagrams $D$. In particular for any two
$x,y\in\mathfrak{K}(F)$ we have the relation
$$\tau (x\star y)=\tau (y)\star \tau (x).$$
For elements in $\mathfrak{C}(F)$ this simplifies to
$$\tau (\alpha \star \beta )=\beta \star \alpha $$
because $\mathfrak{C}(F)$ is invariant under the action of $\tau $.
This is a \textit{hermitian} property of the deformation defined by
the bracket (see \cite{W}).

Note that $g\tau =\tau g$ thus we have a naturally action of 
$\mathcal{M}(F)\times \mathbb{Z}_2\subset \pi_0(Diff(F\times I))$ on $\mathfrak{K}(F)$. 

\vskip .1in

The ring homomorphism $\mathsf{k}[t,t^{-1}]\to \mathsf{k}[[h]]$
defined for rings $\mathsf{k}$ of characteristic $0$ (thus $k\supset
\mathbb{Q}$ and $e^h=\sum_{j=0}^{\infty}\frac{h^j}{j!}$ is defined)
is equivariant with respect to the involution on
$\mathsf{k}[t,t^{-1}]$ and the involution $h\mapsto -h$ on
$\mathsf{k}[[h]]$.  Note also that for each $j$, $<\ >_j$ is
invariant with respect to $\mathcal{M}(F)$ meaning that
$<gD>_j=g<D>_j$, where the action of $\mathcal{M}(F)$ on the right
hand side is the classical action of $\mathcal{M}(F)$ on curve
systems.

\vskip .1in

\noindent \textbf{Remark.} 
It follows immediately from the state sum definition that 
we have the symmetries
$$<\tau (D)>_k=(-1)^k<D>$$
for all $k\geq 0$.

\vskip .5in

\section{The Kauffman resolution algebra.}

For each finite set $S$ let $|S|$ be the number of elements of $S$. 

We will consider connected compact oriented 
surfaces $F$ with $r\geq 0$ boundary components 
equipped with a fixed oriented diffeomorphism (parametrization)
$$\bigcup_rS^1\rightarrow \partial F.$$
This is briefly called a \textit{surface}.
The image of the $i$-th $S^1$ is denoted $\partial_iF$.
The mapping class group $\mathcal{M}(F)$ is the group of 
isotopy classes of diffeomorphisms of $F$ fixing $\partial F$ pointwise.

\vskip .1in

A diagram on $F$ is a pair $(D,C)$ with $D$ a diagram of a regularly immersed 
proper $1$-manifold in $F$ with the usual under-over crossing
information at each crossing, $C\subset F$ is a subset of the the 
the set of crossings of $D$. If the $i$-th boundary component $\partial_iF$
contains $j_i$ boundary points of $D$ then we assume that 
$$D\cap \partial_iF=\bigcup_{\ell=1}^{j_i} \{ e^ {\frac{2\pi
\sqrt{-1}\ell}{j_i}} \}$$

Let $(D,C)$ be isotopic to 
$(D',C')$ if there is an isotopy of $F$ fixing $\partial F$ pointwise
and mapping $D$ to $D'$ and $C$ to $C'$.
The set of isotopy classes of diagrams is denoted $\mathcal{D}(F)$
and we will denote the isotopy class of a pair by $[D,C]$.

\vskip .1in

The set $\mathcal{D}(F)$ naturally decomposes according to $|C|$ and by the 
number of boundary points contained in each of the $r$ components.
Thus
$$\mathcal{D}(F)=\bigcup_{k\geq 0,j\geq 0}\mathcal{D}(F)[k]\{j \},$$
where $k$ is the number elements in $C$ and $j=(j_1,\ldots ,j_r)$ is
a multi-index with $j_i$ the number of boundary points of $D$ in
$\partial_iF$.  Thus $\frac{1}{2}(j_1+\ldots +j_r)$ is the number of
arc components of the diagram $D$.  As usual $j\geq 0$ means
$j_i\geq 0$ for all $i$.  There is also a unique isotopy class of
diagram $\emptyset \in \mathcal{D}(F)[0]\{0 \}$ with $0=(0,\ldots,0)$.

\vskip .1in

For $k\geq 0$ let $\mathcal{D}(F)[k]:=\cup_j\mathcal{D}(F)[k]\{j
\}$. Similar notation applies to the other grading.
Then $\mathfrak{D}(F)$ from section 1 can naturally be
identified with 
the subset of those $(D,C)\in \mathcal{D}(F)[k]\{0 \}$
with $k=|C|$. This set is naturally contained in the subset
$\mathfrak{D}^a(F)$ of all diagrams $[D,C]$ on $F$ for which $C$ is the 
set of \textit{all crossings}. 
We call the elements of $\mathfrak{D}^a(F)$ \textit{real diagrams}.
(Representatives of elements of $\mathfrak{D}^a(F)$ 
possibly contain arc components.)     

\vskip .1in

The above decompositions naturally define the structure of a bi-graded
(with the second grading a multi-grading itself) on the free
$\mathsf{k}$-module with basis $\mathcal{D}(F)$.
Also there is defined the graded submodule 
$\mathsf{k}\mathfrak{D}(F)\subset \mathsf{k}\mathfrak{D}^a(F)$
spanned by isotopy classes of $k$-crossing diagrams in degree
$k\geq 0$. 

\vskip .1in

\noindent \textbf{Definition and Remarks.} We define the 
\textit{Kauffman bracket
module} $\mathfrak{K}^a(F)$ by quotiening the free module 
$\mathsf{k}[t,t^{-1}]\mathfrak{D}^a(F)$ 
by the usual skein relations. A $\mathsf{k}$-module 
basis is given by the set $\mathfrak{C}^a(F)$ of isotopy classes 
of curve systems on $F$ consisting of arbitrary
properly embedded $1$-manifolds with specified boundary conditions as 
described above, but without inessential closed components. 
This set is by definition
$\mathfrak{D}^a(F)[0]\subset \mathcal{D}(F)[0]$.
The right hand side here contains the elements with $C=\emptyset $ but
the diagrams still can have crossings while the left
hand side only contains the real diagrams of this set.
(Note that boundary 
parallel components may be contained in the curve systems.
Thus our notion of curve system is still different from the classical approach.)
We have decompositions:
$$\mathfrak{C}^a(F)=\bigcup_{j}\mathfrak{C}^a(F)\{j \}$$
by specifying the boundary pattern. Also note that
$\mathfrak{C}(F)=\mathfrak{C}^a(F)\{0 \}$
gives the set of curve systems only containing closed 
components considered in section 1. 
The Kauffman bracket module now also  decomposes according to the 
grading. Our version of Kauffman bracket module is a generalization
of the relative Kauffman bracket module \cite{P}.
We like to mention that the weak product defined above actually 
lifts to define a graded product structure on the Kauffman bracket
module:
$$\mathfrak{K}^a(F)\{j \}\otimes \mathfrak{K}^a(F)\{j' \}\rightarrow
\mathfrak{K}^a(F)\{j+j' \}$$ 
by placing a diagram $D$ above a diagram $D'$. Because of isotopy
invariance this is a well-defined product even in this relative case.

\vskip .1in

\noindent \textbf{Question.} Find the
$SL(2,\mathbb{\mathsf{k}})$-interpretation of the relative
Kauffman bracket algebra?  

\vskip .1in

The set-up described above suggests to assign to a curve system 
in $\mathfrak{C}^a(F)\{j \}$ with $\partial F\neq \emptyset$
a \textit{regular} map on the space 
$Flat(F)$ of flat connection on a trivialized
$SL(2,\mathsf{k})$-bundle over $F$
with values in $\mathsf{k}\times SL(2,\mathsf{k})^{|j|}$
and $|j|:=\frac{1}{2}(j_1+ \ldots j_r)$ the number of arcs.
The map is defined using the explicit boundary parametrizations,
basepoint and orientation on each $S^1$ to give an ordering of the set
of arc components of a diagram. In fact, given a flat connection,
each arc components associates the holonomy along the arc.
This defines the mapping to $SL(2,\mathsf{k})^{|j|}$. The unordered
collection of closed components defines a product of functions given
by calculating the traces of holonomies.   
 The bracket relations induces an equivalence relation on this set
 of maps. Details will be discussed in \cite{K}.

\vskip .1in

The mapping class group acts on diagrams and on the free
$\mathsf{k}$-module spanned by diagrams preserving all gradings.
Let $\phi_*$ be the $\mathsf{k}$-endomorphism of degree $(0,0)$ 
induced by $\phi \in \mathcal{M}(F)$.

\vskip .1in

Let $F$ respectively $F'$ be surfaces with $r\geq 1$ respectively
$r'\geq 1$ boundary components. Then we can define a new surface
$F\cup F'$ by glueing the last boundary circles. The boundary
parametrizations can be combined in some obvious way.  Note that all
isotopies and diffeomorphisms fix the boundary and can therefore be
matched. We call $F\cup F'$ the \textit{glueing of $F$ and $F'$}.
Note that $F\cup F'$ has $r+r'-2$ boundary components.  The glueing
operation obviously induces glueing operations of diagrams in the
following way:
\begin{align*}
 \mathcal{D}(F)[k]\{(j_1,\ldots ,j_{r-1},\ell ) \}\times
\mathcal{D}(F_2)[k']\{(j'_1,\ldots ,j'_{r-1}, \ell) \}& \\ 
\to
\mathcal{D}[k+k'](F\cup F')\{(j_1,\ldots ,j_{r-1},j'_1,\ldots,j'_{r-1} \}
\end{align*}
We let $(D,C)\cup (D',C')$ denote the result of glueing 
two diagrams. This is only defined when the number of 
boundary points in the last boundary components match.
This operation is compatible with isotopy and defines
$[D,C]\cup [D',C']$ for given $[D,C]\in \mathcal{D}(F)$,
$ [D',C']\in \mathcal{D}(F')$ which are matching.
The glueing operation extends linearly 
$$\bigoplus_{\ell \geq 0}\mathsf{k}\mathcal{D}(F)\{(j,\ell \}\otimes
\mathsf{k}\mathcal{D}(F')\{(j',\ell \}\to
\mathsf{k}\mathcal{D}(F\cup F')\{(j,j')\}$$
using obvious multi-index notation.

\vskip .1in

We say that a diagram $(D,C)$ or its isotopy class $[D,C]$ is
a \textit{weak product} of diagrams $(D_1,C_1)$ and $(D_2,C_2)$
on a surface $F$ 
if $D$ is given by superimposing $D_1$ and $D_2$ to form a diagram
$D_1\triangleright D_2$ with all crossings of the form $D_1$ over $D_2$.
We will have $C=C_1\cup C_2$ thus ignore all new crossings 
of $D_1$ with $D_2$.
In general this requires modifications in the boundary using 
natural diffeomorphisms of the circle to be able to take the union 
in the boundary. Thus actually $(D_i,C_i)$ will be modified by 
isotopy of $F$ in a neighbourhood of the boundary circles. We 
will not formalize this construction at this point. 
In general superimposing diagrams is not a 
well-defined operation on isotopy classes of diagrams. Thus  
$[D_1,C_1]\triangleright [D_2,C_2]$ will be the notation for \textit{any}
weak product resulting from representatives $(D_1,C_1)$ and
$(D_2,C_2)$. Similarly there are defined \textit{strong products} 
by adding all crossings of $D_1$ with $D_2$ to the set of 
crossings $C$ of $D$.
 
\vskip .1in

Now let $\mathcal{E}(F)[i]$ denote the set of
$\mathsf{k}$-endomorphisms of $\mathsf{k}\mathcal{D}(F)$ of bi-degree
$(-i,0)$, i.e. $\rho \in \mathcal{E}(F)[i]$ if 
$\rho (\mathcal{D}(F)[k]\{j \})\subset \mathcal{D}(F)[k-i]\{j \}$
for all $k\geq 0$ and all $j\geq 0$. 
(Here we set $\mathcal{D}(F)[k]=0$ for $k<0$.)
Then define the graded algebra 
$$\mathcal{E}(F):=\bigoplus_{i\geq 0}\mathcal{E}(F)[i].$$
This is a subalgebra of the $\mathsf{k}$-algebra of all 
$\mathsf{k}$-endomorphisms of the $\mathsf{k}$-module
$\mathsf{k}\mathcal{D}(F)$.

\vskip .1in

We will call the assignment of a graded subalgebra 
$$F\mapsto \mathcal{R}(F)\subset \mathcal{E}(F)$$
a \textit{Kauffman resolution functor} if it satisfies the following 
conditions \textbf{1. - 6.}
(This is formally a functor if we consider the category of surfaces as objects 
and morphisms between surfaces compatible with boundary 
parametrizations. In fact, morphisms will induce
$\mathsf{k}$-homomorphisms
of diagram algebras and thus homomorphisms of their graded 
endomorphism algebras in a natural way.)

\vskip .1in

\noindent \textbf{1.\ Mapping class invariance.} Given $\rho \in
\mathcal{R}(F)$ and $\phi \in \mathcal{M}(F)$ then 
$$\phi_* \circ \rho =\rho \circ \phi_*$$  

\vskip .1in

\noindent \textbf{2.\ Glueing.} Let surfaces $F_1,F_2$ be given (both 
with nonempty boundary) such that the glueing surface $F_1\cup F_2=F$
is defined. 
Then there 
exists a unique restriction homomorphism of degree $0$:
$$\mathcal{R}(F)\ni \rho \mapsto (\rho |F_1) \in \mathcal{R}(F_1)$$
such that for all diagrams $[D,C]$ with $C\subset F_1$
such that $D$ intersects the image of the distinguished circle
in $F$ transversely:
$$\rho [D,C]=(\rho |F_1)\left([D\cap F_1,C]\right)\cup [D\cap F_2,\emptyset].$$

\vskip .1in

\noindent \textbf{3.\ Weak product.} For any two diagrams 
$[D',\emptyset], [D,C]\in \mathcal{D}(F)$ for which weak products 
are defined we have:
$$\rho ([D',\emptyset]\triangleright [D,C])=[D',\emptyset]
\triangleright \rho [D,C].$$
Here the right hand side is interpreted as the linear combination
of weak products (be aware that this operation is not a well-defined
operation on isotopy classes) of some representative $(D',\emptyset)$
with the representatives of the terms in $\rho [D,C]$. 
The same identities are supposed to
hold with the order of $D',D$ switched.

\vskip .1in
  
Note that for $\rho \in \mathcal{R}(F)[i]$ and $[D,T]\in
\mathcal{D}(F)[i]$ we have $\rho [D,T]\in
\mathsf{k}\mathcal{D}(F)[0]$, 
so we write formally:
$$\rho [D,T]=[\overline{\rho }[D,T],\emptyset],$$
because the crossing information is empty for all terms.
In fact, in general let $D_i$ be a sequence of diagrams with 
the same set $C\subset F$ of crossings with respect to a choice 
of representative diagrams. Then
$$\sum \lambda_i [D_i,C]=[\sum \lambda_iD_i,C]$$
has a well-defined meaning. 
\vskip .1in

\noindent \textbf{4.\ Generalized divergence.} This property results
from the idea that pairs $(D,C)$ can be interpreted formally with
$D$ a function and $C$ a set of variables of the function.
Let $F$ be a surface and $[D,C]\in \mathcal{D}(F)$,
$\rho \in \mathcal{R}(F)[i]$. Then
$$\rho [D,C]=\sum_{T\subset C,|T|=i}[\overline{\rho
}[D,T],C\setminus T].$$
In order to make sense of the right hand expression we need to justify
that $C\setminus T$ is \textit{naturally} a subset of all the diagrams
in $\rho [D,T]\in \mathsf{k}\mathcal{D}(F)$.  
This is the essential technical step and follows from the glueing axiom \textbf{2.}
above applied to a splitting of $F$ along a curve separating the 
crossings in $T$ from the crossings in $C\setminus T$. 
Here we split the diagram into two diagrams $D_1\subset F_1$ with
all the crossings of $T$ in $D_1$, and $D_2\subset F_2$ containing the
crossings of $C\setminus T$.
Then the glueing axiom implies that complete diagram $D_2$ is glued 
back to the terms in $(\rho |F_1)[D_1,T]$ thus the result naturally
contains $C\setminus T$.

\vskip .1in

\noindent \textbf{5.\ Skein relation.} For each surface $F$ there 
exist two module epimomorphisms of degree $-1$:
$$\tilde{\zeta},\tilde{\iota} \co \mathcal{R}(F)\rightarrow
\mathcal{R}(F)$$
satisfying
$$\tilde{\zeta}\tilde{\iota }=\tilde{\iota }\tilde{\zeta },$$
and such that $\tilde{\zeta}(\rho), \tilde{\iota}(\rho)$ are linearly independent 
for all $\rho\in \mathcal{R}(F)$.
The homomorphisms $\tilde{\zeta}, \tilde{\iota }$ have to satisfy 
that for all $\rho \in \mathcal{R}(F)[i]$ and 
$[D,C]\in \mathcal{D}(F)[i]$ the following skein relation holds:
$$\rho [D_+,C]+(\tilde{\zeta }(\rho ))[D_0,C\setminus \{+ \}]+
(\tilde{\iota}(\rho))
[D_{\infty},C\setminus \{+ \}]=0\in \mathsf{k}\mathcal{D}(F),$$
where $+$ is a crossing of $D$ in $C$, and $D_0, D_{\infty}$ are 
the usual resolutions.  

\vskip .1in

\noindent \textbf{6.\ Vacuum condition.} For all $F$ and $\rho \in
\mathcal{R}[0]$ there exists a constant $\theta \in \mathsf{k}$ such that
$$\rho ([\emptyset ])=\theta [ \emptyset ].$$
Moreover, for each $\theta \in \mathsf{k}$ there exists some
$\rho \in \mathcal{R}(F)[0]$ with this property. 

\vskip .1in

After this long technical preparation we can now state the main result
of this section. 

\vskip .1in

\noindent \textbf{Theorem 1.} \textit{
There exists at most one assignment 
$F\mapsto \mathcal{R}(F)\subset \mathcal{E}(F)$  
satisfying properties} \textbf{1.~-~6.} \textit{above. 
Moreover, under the assumption that the functor exists, the 
value of $\rho \in
\mathcal{R}(F)[i]$ on a diagram $[D,C]$ with $|C|=i$ is determined 
by a state summation over Kauffman states on $C$ with the 
coefficients in $\mathsf{k}$ determined only by the number of
$\infty$-states of a state.}

\vskip .1in

\noindent \textit{Proof.} It follows from the 
the weak product property for
$(D,C)=\emptyset $ that 
$$\rho [D',\emptyset]=[D',\emptyset]\cdot \rho
[\emptyset ].$$ Because of the grading, $\rho $ can possibly be
nonzero on $[\emptyset]$ only for $\rho \in
\mathcal{R}[0]$. Otherwise $\rho [\emptyset]=0$ and thus $\rho $
vanishes on $\mathcal{D}(F)[0]$.

Now in general $\rho \in \mathcal{R}(F)[i]$ vanishes on
$\mathcal{D}(F)[j]$ for $j<i$ because of the grading.  Moreover, the
divergence property determines $\rho $ on $[D,C]$ with $|C|=j>i$
from the values on diagrams in $\mathcal{D}(F)[i]$.  The result then
is proved by induction over $i$ using the skein relation in
combination with the glueing property.  More precisely it follows
from the skein relation that $\rho $ is determined by
$\tilde{\zeta}(P)$ and $\tilde{\iota}(P)$.  For $i=1$ it follows
that $\tilde{\zeta}(P)$, respectively $\tilde{\iota }(P)$, acts on
$\mathcal{D}(F)[0]$ as multiplication by a constant in
$\mathsf{k}$. Thus due to the linear independence the degree $1$
resolution is determined by two numbers $a_0,a_1\in \mathsf{k}$.  Of
course in this case the coefficients of $\rho [D,C]$ are determined
by the number $0$, respectively $1$, of $\infty $-states. For the
induction step from $i-1$ to $i$ we first note that by induction
hypothesis $\tilde{\zeta}(\rho)$ is determined by $i-1$ numbers
giving the contribution of a state sum with $i$ $\infty$-markers in
a state summation over $D_0$.  We will prove that
$\tilde{\iota}(\rho)$ is determined by just one more
coefficient. Now $\tilde{\iota }(\rho )$ is determined also by $i-1$
numbers, and in fact from
$\tilde{\zeta}\tilde{\iota}(\rho)=\tilde{\iota}\tilde{\zeta}(\rho)$
and $\tilde{\iota}\tilde{\iota}(\rho)$. Now the first contribution
is already known from $\tilde{\zeta }(\rho )$.  We can iterate the
application of $\tilde{\iota }$ and $\tilde{\zeta}$ and use
induction hypothesis to reduce to $\tilde{\iota}^{i}(\rho )$, which
is of degree $0$ and thus determined by a single coefficient.  Since
in the applications of $\tilde{\iota }$ we smooth a crossing each
time it is obvious that this coefficient is determined by the number
of $\infty$-states. $\square$

It is the main result of the next section that the algebra
$\mathcal{R}(F)$ exists for each compact connected oriented surface
and is naturally graded isomorphic to the the polynomial 
algebra $\mathsf{k}[z,w]$. 

\vskip .1in

\noindent \textbf{Remarks.} (a) It does not seem to be possible to
characterize the the resolution algebras $\mathcal{R}(F)$ without
extending the axiomatic to surfaces with boundaries possibly
containing arc components, even if we are finally only interested in
closed surfaces.  The crucial property is the glueing property which
defines the locality of the operations.  The glueing axiom is
necessary just to formulate the crucial divergence property which
\textit{localizes} the operation of $\rho$ of degree $i$ on
$i$--crossing diagrams. Similarly it seems difficult to develop the
axiomatic characterization without the flexibility in the grading by
crossing numbers.

\noindent (b) Suppose we consider the case $j=0$, i.e. no arcs on a
closed surface $F$.  In this case the vacuum condition can be
actually deduced from the other properties.  Then we know that $\rho
[\emptyset]$ is a $\mathsf{k}$-linear combination of elements of
$\mathcal{D}(F)[0]\{0\}$. Because $F$ is closed it follows easily
from naturality applied to some Dehn twist $\phi $ of sufficient
high order that only the empty diagram can appear in the linear
combination.  In fact $\phi_*[\emptyset] = [\emptyset]$ while
$\phi$ can be chosen such that the finite linear combination of
nonempty curve systems is \textit{not} fixed under $\phi_*$.

\vskip .1in

The involution $\tau$ given by changing crossings obviously extends
to an involution of $\mathcal{D}(F)$ of degree $(0,0)$ by changing the
crossings of $C$ but fixing the other crossings.

\vskip .5in

\section{Existence of the resolution algebra.}

We will consider the sequence of brackets for $j\geq 1$, also denoted
$$<\ >_j\co\mathsf{k}\mathfrak{D}^a(F)\rightarrow \mathsf{k}\mathfrak{C}^a(F),$$
defining the Kauffman bracket as in section 1. But we now work
in the more general case of diagrams and skein modules possibly
containing proper arcs.
 
\vskip .1in

In the following the grading of the the polynomial algebra
$\mathsf{k}[z, w]$ is given by the total degree. We like to point
out that the variables $z,w$ correspond to the state maps $\iota,
\zeta$ and its associated operator versions $\tilde{\iota },
\tilde{\zeta}$ from section 3.

\vskip .1in

\noindent \textbf{Theorem 2} \textit{For each surface $F$ there 
is a graded homomorphism of $\mathsf{k}$-algebras
$$\chi : \mathsf{k}[z,w]\rightarrow \mathcal{E}(F)$$
If $\mathsf{k}$ is a ring of characteristic $0$ then $\chi $ is
injective and the image is a Kauffman resolution algebra of $F$.}

\vskip .1in

\noindent \textit{Proof.}  First let 
$$\mathfrak{c}\co \mathsf{k}[z,w]\rightarrow \mathsf{k}[z,w]$$
be the algebra homomorphism defined by mapping $z$ to $zw $ and $w$
to $w$.  The image of $\mathfrak{c}$ is the algebra of polynomials
in $\mathsf{k}[z,w]$ of the form $\sum p_i(z)w ^i$ with
$deg(p_i)\leq i$. Note that the image of a homogeneous polynomial $P$
of degree $k$ is of the form $p(z)w^k$ with a polynomial $p$ in $z$
of degree $\leq k$.  We will define $\chi_{P}$ in terms of the
homogeneous components of $\mathfrak{c}(P(z,w))$. Given $k$ we will
first define $\chi_k$ on polynomials $p(z)$ of degree
$\leq k$.  Note that $p(z)=a_o+a_1z+\ldots +a_kz^k$ is determined by
the mapping
$$\mathfrak{p}\co \{0,1,\ldots ,k\}\rightarrow \mathsf{k}$$
with $\mathfrak{p}(j)=a_j$.  

Now let $(D,C)\in \mathcal{D}(F)$.  A $k$-state on $C$ is a choice of
$k$-element subset $T$ of $C$ and a mapping 
$\sigma\co T \rightarrow \{0,\infty\}$.  We denote the set of all
$k$-states on $(D,C)$ by $\Sigma_k(C)$. 
Then define $\chi_k\co \mathsf{k}[z] \to \mathcal{E}(F)$ by
$$
\chi_k(p)[D,C]:=(-1)^k \sum_{\sigma\in\Sigma_k(C)}
     \mathfrak{p}(\iota(\sigma))[D(\sigma ),C(\sigma)]
$$
where $D(\sigma )$ is the diagram which results from $D$ by
smoothing the crossings in the domain $T$ of $\sigma $ as determined
by $\sigma$, and $C(\sigma )=C\setminus T$ for each state $\sigma\co
T \rightarrow \{0,\infty \}$.  For each natural number $k$ let
$\pi_k\co \mathsf{k}[z,w]\to \mathsf{k}[z]$ be the map that sends a
polynomial $P(z,w)$, considered as an element of $\mathsf{k}[z][w]$,
to its $k$\,th coefficient.  Now for a general polynomial
$P(z,w)\in\mathsf{k}[z,w]$ define $\chi$ by
$$
\chi = \sum_{j=0}^k\chi_j\circ \pi_j \circ \mathfrak{c}
$$

Consider two homogeneous polynomials $P,Q\in \mathsf{k}[z,w]$ of degree $j$,
respectively $k$. Let $\mathfrak{p}$, respectively $\mathfrak{q}$,
be the corresponding function $\{0,1,\ldots ,k\}\rightarrow\mathsf{k}$.
Let $\mathfrak{r}\co \{0,1,\ldots ,j+k\}\rightarrow \mathsf{k}$ be the
function determined by the polynomial $PQ$. Note that if $\mathfrak{c}(P)=p(z)w^j$ and 
$\mathfrak{c}(Q)=q(z)w^k$ with $deg(p(z)\leq j$ and $deg(q(z))\leq k$ then 
$$\mathfrak{c}(PQ)=\mathfrak{c}(P)\mathfrak{c}(Q)=p(z)q(z)w^{j+k}.$$ 
Thus
$$\mathfrak{r}(i)=\sum_{i=i_1+i_2,i_1\leq j,i_2\leq k}\mathfrak{p}(i_1)\mathfrak{q}(i_2).$$
It follows from the definition that
$$
  \chi_P\chi_Q[D,C] = (-1)^{j+k}\sum_{\sigma \in \Sigma_j(C(\tau))}
 \mathfrak{p}(\iota(\sigma))
  \sum_{\tau \in \Sigma_k(C)}
  \mathfrak{q}(\iota(\tau)) [(D(\tau ))(\sigma),C((\tau )(\sigma ))]
$$
which is equal to the state sum
$$(-1)^{j+k} \sum_{\eta \in \Sigma_{j+k}(C(\tau))} 
\mathfrak{r}(\iota(\eta)) [D(\eta),C(\eta )].$$ 

It is obvious from the definition that $\chi_P$ is compatible with
the action of $\mathcal{M}$.  It is also clear that the image
satisfies \textbf{1.\ - 6.} of section 3.  It remains to show
injectivity for rings $\mathsf{k}$ of characteristic $0$. We only
have to show that $P\neq 0$ implies that $\chi_P\in \mathcal{E}(F)$
is not the trivial endomorphism.  Let
$$P=P_i+P_{i+1}+\ldots +P_N$$
be the decomposition of $P$ into homogeneous components with
$P_i\neq 0$.
Consider the following $i$-crossing diagram of a circle on a disk in
$F$. 

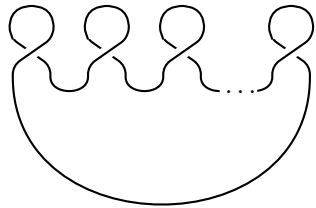
\begin{center}
\unitlength 0.5cm
\psset{unit=\unitlength}
  \begin{pspicture}(-.5,-3.5)(9,2)
    \psline[linearc=.45](1,0)(1,.3)(0,1)(-.1,1.3)
    \psline[linearc=.45,border=.15](0,0)(0,.3)(1,1)(1.1,1.3)(1,1.8)(.5,1.9)(0,1.8)
    (-.1,1.3)(0,1) \psbezier
    (1,0)(1,-.5)(2,-.5)(2,0)
    \psline[linearc=.45](3,0)(3,.3)(2,1)(1.9,1.3)
    \psline[linearc=.45,border=.15](2,0)(2,.3)(3,1)(3.1,1.3)(3,1.8)(2.5,1.9)(2,1.8)
    (1.9,1.3)(2,1) \psbezier
    (3,0)(3,-.5)(4,-.5)(4,0)
    \psline[linearc=.45](5,0)(5,.3)(4,1)(3.9,1.3)
    \psline[linearc=.45,border=.15](4,0)(4,.3)(5,1)(5.1,1.3)(5,1.8)(4.5,1.9)(4,1.8)
    (3.9,1.3)(4,1) \psbezier
    (5,0)(5,-.5)(6,-.5)(6,0) \psbezier
    (5.9,0)(5.9,-.5)(6.9,-.5)(6.9,0)
    \psframe[linecolor=white,fillcolor=white,fillstyle=solid](5.5,-1)(6.5,1)
    \rput(6,-.4){ \ldots}
    \psline[linearc=.45](7.9,0)(7.9,.3)(6.9,1)(6.8,1.3)
    \psline[linearc=.45,border=.15](6.9,0)(6.9,.3)(7.9,1)(8.0,1.3)(7.9,1.8)(7.4,1.9)
    (6.9,1.8)(6.8,1.3)(6.9,1) 
    \psbezier
    (0,0)(0,-4.5)(7.9,-4.55)(7.9,0)
  \end{pspicture}
\end{center}

We let $(D,C)$ be this diagram with $C$ the set of all crossings of
$D$. Then $\chi (P_j)$ is trivial on $\mathcal{D}(F)[i]$ for $j>i$
because there are no $j$-element subsets of the set of crossings.
The smoothing according to some $i$-state with $\sigma^{-1}(\infty
)=\ell$ is an $\ell+1$-component diagram in the disk in $F$. There
are precisely ${i \choose \ell}$ states of this form, and all give
rise to the same diagram.  Thus each nontrivial coefficient
$a_{\ell}$ in $P$ will contribute a coefficient ${i \choose \ell}
a_{\ell}\neq 0$ in $\chi_P[D,C]$, which does not cancel with any other
contribution.  $\square$ 

\vskip .2in

In the following we only consider $(D,C)$ with $C$ the set of all
crossings of $D$. In this case we only write $D$ both for a
representative diagram and its isotopy class.

\vskip .1in

\noindent \textbf{Examples.} (a) For a constant polynomial $P=a_0\in
\mathsf{k}\subset \mathsf{k}[z,\zeta ]$ 
the sum is over the single $0$-element subset of the set of
crossings of $D$ 
and contributes $a_0D$ because
$|\sigma ^{-1}(\infty )|=0$. Thus $\chi_{a_0}$ is just multiplication by $a_0$. 
\noindent (b) For $k=1$ and $P=w-z$ thus $\mathfrak{c}(P)=(1-z)w$
we know that $\mathfrak{p}(0)=1$ and $\mathfrak{p}(1)=-1$. 
Thus
$$\chi_P(D)=-\sum_{\textrm{crossings} \ p \ \textrm{of} \ D }(D_{p,0}-D_{p,\infty}),$$
which is just the Poisson bracket defined in section 1. 

\noindent (c) If $D$ is a diagram with $k$ crossings and $P$ is a
homogeneous polynomial of degree $>k$ then $\chi_P(D)=0$.  Thus if
diagram resolutions are considered to operate like differential
operators on functions, this behaviour very much suggests
$k$-crossing diagrams to correspond to polynomial functions of
degree $k$. We say that each element in $\mathcal{R}(F)$ has
\textit{finite support} (with respect to the $[\ ]$-grading
respectively number of crossings if restricted to real diagrams.)

\noindent (d) We can apply theorem 2 to the ring
$\mathsf{k}[t,t^{-1}]$ itself. Define for $k\geq 0$ the sequence of
polynomials
$$P[k](z,w)=t^kw^k+t^{k-2}w^{k-1}z+\ldots +t^{-k}z^k\in \mathsf{k}[t,t^{-1}][z,w].$$
Then for each diagram $D$ with $k$ crossings
$$\chi_{P[k]}(D)=<D>'\in \mathsf{k}[t,t^{-1}]\mathfrak{D}^a(F)[0].$$
Then there is a natural homomorphism
$$\mathsf{k}[t,t^{-1}]\mathfrak{D}^a(F)[0]\rightarrow \mathsf{k}[t,t^{-1}]\mathfrak{C}^a(F)$$
defined by mapping the curve system $\gamma $ to
$(-t^2-t^{-2})^{\mu(\gamma )}\gamma_0$, where $\mu $ is the number
of trivial components in $\gamma $ and $\gamma_0$ is the result of
discarding the trivial components from $\gamma $.  The definition of
$<D>$ thus is separated into two steps. Similarly we will separate
the calculation of $<D>_j$ for all $j\geq 0$.  The Kauffman bracket
respectively its extension can be considered as an operator with
\textit{infinite} support:
$$\mathcal{D}(F)\rightarrow \mathcal{D}(F)[0]$$
of the form
$$\hat{P}=\sum_{k=0}^{\infty}\chi_{P[k]}\circ \Pi_k\in \mathcal{E}(F)$$
where 
$$\Pi_k\co \mathcal{D}(F)\rightarrow \mathcal{D}(F)[k]$$
is the projection onto the $k$-th grading module. It maps
$\mathcal{D}(F)[k]\subset \mathcal{D}(F)[k]$ by the identity 
and mapping all other $\mathcal{D}(F)[j]$ trivially.
Thus $\Pi_k\in \mathcal{E}(F)[0]$ but $\Pi_k\notin \mathcal{R}(F)$.
Theorem 2 can be considered as an \textit{operator expansion} of the 
$\mathsf{k}[t,t^{-1}]$-operator with infinite support in terms of 
finite support $\mathsf{k}$-operators.

\vskip .2in

The following result is immediate from the definition of $\chi$
and generalizes the skein relation from section 3. 

\vskip .1in

\noindent \textbf{Theorem 3.} \textit{For each homogeneous
  polynomial $P$ and diagram $D$ with usual Kauffman triple $D_+,
  D_0, D_{\infty}$, the following $SL(2,\mathbb{C})$-skein relation
  holds:}
$$
\chi_P(D_+-D_0-D_{\infty})+\chi_{(P-a_kz^k)w^{-1}}(D_0)+\chi_{(P-a_0w^k)z^{-1}}(D_{\infty})=0.
$$  

\vskip .1in

In order to be able to find combinatorial expressions of
$$\mathfrak{D}^a(F)\ni D\rightarrow <D>_j \in \mathsf{k}\mathfrak{C}^a(F)$$
in terms of our algebra $\mathcal{R}(F)$, we need to define certain
projection homomorphisms into $\mathsf{k}\mathfrak{C}^a(F)$.

\vskip .1in

Let $\varphi \co \mathbb{N}=\{0,1,2,\ldots \}\rightarrow \mathsf{k}$ be any map. 
Define 
$$\varphi_* \co \mathsf{k}\mathcal{D}(F)\rightarrow \mathsf{k}\mathcal{D}(F)$$
by
$$\varphi_*[D,C]=(-1)^{|C|}\sum_{\textrm{states} \ \sigma \ \textrm{on} \
  D} \varphi(\mu(\sigma ))[D(\sigma ),\emptyset ],$$ where as before
$|C|$ is the number of crossings of $D$, $\mu (\sigma )$ is the number
of trivial components of the smoothing of $D$ using $\sigma $, and
$D(\sigma )$ is the diagram resulting from discarding the trivial components from
this smoothing.

\vskip .1in

If applied to $(D,C)\in \mathfrak{D}^a(F)$ then
$$\varphi_*[D,C]\in \mathsf{k}\mathfrak{C}^a(F).$$

\vskip .1in

For the function $\varphi (i)=(-2)^i$ we have
$\phi_*(D)=<D>_0$.
 
Also note that $\varphi_*(P(D))=<D>_1$ if $P=w-z$.   

\vskip .1in

Suppose that $\mathsf{k}$ has characteristic $0$.  Our basic
projections are defined from the sequence of maps
$$\varphi _j \co \mathbb{N}\rightarrow \mathbb{Q}\subset \mathsf{k}$$
defined by $\varphi_j(i)$ is the the coefficient of $h^{2j}$ in the expansion 
of $(-t^2-t^{-2})^i$ with $t=e^h$. 
Then
$$
  \varphi_j(i) =(-1)^i\frac{2^{2j}}{(2j)!}\sum_{k=0}^i{i \choose k}(2k-i)^{2j}
$$

In particular
$$\varphi_0(i)=(-1)^i\sum_{k=0}^i{i\choose k}=(-2)^i,$$
and
$$\varphi_1(i)=-(-2)^{i+1}i.$$

\vskip .1in

\noindent \textbf{Remark.} It is important to observe that
$\varphi_*[D,C]$ does \textit{not} depend on the over-undercrossing information
of the crossings in $C$ for any map $\varphi $.

\vskip .5in

\section{Combinatorial expansion.}

We are now ready to state the main result of the paper.
We assume that $\mathsf{k}$ is a ring of characteristic zero and 
identify $\mathbb{Q}\subset \mathsf{k}$.

\vskip .1in

\noindent \textbf{Theorem (technical version).} \textit{For each $k\geq 0$
there exists a polynomial 
$$P_k\in \mathbb{Q}[z,w]$$
of degree $k$ (but not homogeneous for $k>1$) such that
$$<D>_k=\sum_{j=0}^{\lfloor \frac{k}{2} \rfloor}(\varphi_j)_*\circ
\chi (P_{k-2j})(D)\in \mathbb{Q}\mathfrak{C}^a(F)\subset \mathsf{k}\mathfrak{C}^a(F)$$
for all $[D]\in \mathfrak{D}^a(F)$. Moreover the homogeneous degree
$k$ component
of $P_k$ is given by
$$w^k-zw^{k-1}+z^2w^{k-2}-\ldots +(-1)^kz^k.$$}

\vskip .1in

The terms with $j>0$ are the loop correction terms. They play an
important tole in section 6.

\vskip .1in

\noindent \textbf{Corollary.} \textit{Let $\alpha, \beta $ be two simple
essential loops on $F$ and $\alpha \cdot \beta $ be a diagram of
$\alpha $ over $\beta $. Then for all $k\geq 0$:
$$\lambda_k(\alpha ,\beta )=\sum_{j=0}^{\lfloor \frac{k}{2}\rfloor}
(\varphi_j)_*\chi (P_{k-2j})(\alpha \cdot \beta ).$$}

\vskip .1in

The corollary also holds more generally for products for $\alpha,
\beta $ possibly proper arcs.

\vskip .1in

The polynomials $P_k$ are given for small $k$ by

$$P_0=1$$

$$P_1=w-z$$

$$P_2=(w^2-zw+z^2)+\frac{1}{2}(w+z)$$

$$P_3=(w^3-zw^2+z^2w-z^3)+(w^2-z^2)-\frac{1}{6}(w-z)$$

A polynomial $P(z,w)$ is called \textit{symmetric} respectively 
\textit{anti-symmetric} if $P(w,z)=P(z,w)$ respectively 
$P(w,z)=-P(z,w)$. 

\vskip .1in

\noindent \textbf{Proposition.} \textit{The polynomials $P_k$ are symmetric
for $k$ even and anti-symmetric for $k$ odd.} 

\vskip .1in

\noindent \textit{Proof.}  This follows from the remark at the end of
section 2 together with the obvious fact that if $\bar{P}(z,w)=P(w,z)$
then $\chi_{\bar{P}}(D)=\chi_P(\tau (D))$. $\square$ 
 
\vskip .1in

\noindent \textbf{Remarks.} 
(a) Note that $\mathsf{k}[z,w]$ has a module splitting in symmetric and 
anti-symmetric polynomials. Of course, as an algebra it is generated by 
$z-w$ and $z+w$. While the first polynomial corresponds to the 
Goldman Poisson bracket, the symmetric generator does not define
an algebraic structure on $\mathcal{C}(F)$. Note that the
commutator $[\alpha ,\beta]=\alpha \star \beta - \beta \star \alpha$
expands in terms of only anti-symmetric resolution operations
(i.e. coming from anti-symmetric polynomials via $\chi $).
But this module is not spanned by $z-w$ alone. The important point is  
that the symmetry defines an additional $\mathbb{Z}_2$-grading on the algebra.
Thus our result shows that essentially the coefficients
$\lambda_j$ can all be deduced from first order operations. 
But the existence of higher order operations seems to be related
to the associativity of the $\star $-product. Compare \cite{CS}
and also the recent work of Abouzaid on the Fukaya category of 
higher genus surfaces \cite{A}. 

\noindent (b) The homogeneous component of maximum degree in each order
is reminiscent of the natural star products of Gutt and Rawnsley 
\cite{GR}. 

\vskip .1in

\noindent \textit{Proof of the theorem.}
The main idea of the proof is already contained in \cite{BFK}.
We discuss a state model computation for $(\varphi_j)_*\circ \chi_P$
and $P$ a homogeneous polynomial of degree $k$.
Note that states for the computation here consist of pairs 
consisting of a state on a $k$-element subset and a state on the remaining
set of crossings. The set of those states maps onto the set of 
Kauffman states. In fact, many states for $(\varphi_j)_*\circ \chi_P$
will contribute to the same Kauffman state. Recall that a Kauffman state 
has $\zeta (\sigma )$ $0$-states and $\iota (\sigma )$ $\infty$-states. 
Recall that the polynomial $P$ is determined by the sequence of 
coefficients $a_0,\ldots ,a_k$ giving the weights associated to
states on $k$-element subsets where $a_j$ is the weight corresponding
to a state with $j$ $\infty $-markers. Now there are ${\iota (\sigma
)\choose j}{\zeta (\sigma ) \choose k-j}$ different states, which will
all give rise to the same Kauffman state and will be have weight
$a_j$. 
The idea is to work within a Kauffman state and  expand using
the functions $\zeta $, $\iota $ and $\mu $ on states as variables.
 
In the calculation of $<\ >$ the term of order $h^k$ is calculated 
from the expansions of $e^{h(\zeta -\iota )}$ and the expansion
of $(-e^{2h}-e^{-2h})^{\mu }$ by collecting the terms whose degree
adds up to $k$. We will consider that summand with order $k$ in
$e^{h(\zeta -\iota )}$ and order $0$ in $(-e^{2h}-e^{-2h})^{\mu }$.
Note that this means that the contribution from the trivial components
will give multiplication by $2^{\mu }$ precisely as in the definition
of $(\varphi_0)_*$.
Then it is easy to see that the other terms are calculated from
the polynomials $(\varphi_j)_*\circ P_{k-2j}$.
Note that
$$e^{h(\zeta -\iota )}=\sum_{k=0}^{\infty}\frac{1}{k!}(\zeta -\iota
)^kh^k,$$
so in order $k$ we have to calculate
$$\frac{1}{k!}\sum_{j=0}^k{k\choose j}(-1)^j\zeta ^{k-j}\iota ^j$$
Consider the $j$-th term in this sum with coefficient 
$$c_j=(-1)^j\frac{\iota ^j\zeta ^{k-j}}{j!(k-j)!}.$$
This has to be compared with the term
$${\iota \choose j}{\zeta \choose {k-j}},$$
which is equal to
$$\frac{1}{j!(k-j)!}\frac{\iota !}{(\iota -j)!}\frac{\zeta !}{(\zeta -(k-j))!}$$ 
or
$$\frac{1}{j!(k-j)!}\iota (\iota -1)\ldots (\iota -j+1)\zeta (\zeta
-1)\ldots (\zeta-(k-j)+1).$$
This is in homogeneous order $k$ in $\zeta $ and $\iota $
precisely $(-1)^jc_j$. 
The result now follows by choosing the coefficients of $P_k$ as in
the theorem. Then the highest homogeneous terms coincide and we have 
expanded $<\ >_k$ in terms of the degree $k$-term $P_k^{(k)}$ of
$P_k$ as given above and lower order terms. These lower 
order terms of $P_k$ are necessary to compensate for the additional
contributions of $(\varphi_0)_*\circ P_k^{(k)}$. 
$\square$

\vskip .1in

The proof of the theorem shows that the explicit calculation of the
polynomials while easy in principle, is in fact a tedious exercise
in binomial combinatorics. It should be very interesting to have an
inductive way of calculation which then could be considered as a
combinatorial \textit{Baker-Campell-Hausdorff} expansion, hopefully
related with geometric structures on the representation variety, see
also \cite{Ku}.

\section{Differentiability of the deformation.}

Let $A$ be a $\mathsf{k}$-algebra. Let $D\subset End_{\mathsf{k}}(A)$
be a filtered subalgebra, i.e. a sequence of sub modules
$$D_0\subset D_1\subset \ldots \subset D_p\subset \ldots \subset
End_{\mathsf{k}}(A)$$
such that the restriction of the multiplication of $A$ satisfies
$$D_i\cdot D_j\subset D_{i+j}$$
for all non-negative integers $i,j$.
For $i\geq 0$, elements of $D_i\setminus D_{i-1}$ are called $D$-operators of
order $i$. Let 
$$D:=\bigcup_{i\geq 0}D_i.$$

\vskip .1in

Recall that a $\star $-product on a $\mathsf{k}$-algebra $A$ is a
$\mathsf{k}[[h]]$-bilinear map
$$A[[h]]\otimes A[[h]]\rightarrow A[[h]],$$
thus is determined by 
$$A\otimes A\rightarrow A[[h]],$$
and thus by a sequence of $\mathsf{k}$-bilinear homomorphisms
$$\lambda_k\co A\otimes A\rightarrow A.$$
for $k\geq 0$.

\vskip .1in

\noindent \textbf{Definition.} A 
$\star $-product on $A$ is called \textit{$D$-differentiable}
if for all $k\geq 0$ the restrictions of the corresponding
sequence of $\mathsf{k}$-bilinear 
homomorphisms
$$\lambda_k\co A\otimes A\rightarrow A$$
to each variable are $D$-operators of order $\leq k$. 

\vskip .2in

Note that for $\star $-products with $\lambda_k$
symmetric or anti-symmetric for each $k$, it suffices to consider the
restriction to the second (or first) variable. 

\vskip .2in

The above definition generalizes the usual definition of
differentiability
of deformations of algebras using the filtered subalgebra
$\mathcal{D}=\mathcal{D}(A)$ of differential operators defined as follows, see
\cite{MR}:
Let $\mathcal{D}(A)_0:=A$ acting by multiplication of $A$ on $A$ and 
inductively for $p\geq 1$ 
$$\mathcal{D}(A)_p:=\{f\in End_{\mathsf{k}}(A)|fa-af\in \mathcal{D}(A)_{p-1}
\ \textrm{for \ all} \ a\in A \}$$   
  
\vskip .1in

In our situation we have $A=\mathsf{k}\mathfrak{C}(F)$
equipped with the $\star $-product induced by the Kauffman 
bracket in $F\times I$.
Then following \textit{combinatorial} filtration
is naturally defined in this case. 
We will write $\varphi_i$ for $(\varphi_i)_*$ to simplify notation.
Let
$D_p$ be the set of those $f\in End_{\mathsf{k}}(A)$, that can be
written as
linear combinations of 
homomorphisms
$$\beta \mapsto
\varphi_{r_{\ell}}\chi_{Q_{\ell}}(\alpha_{\ell -1}\triangleright
\ldots (\alpha_2\triangleright
\varphi_{r_1}\chi_{Q_1}(\alpha_0\triangleright
\varphi_0\chi_{Q_0}\beta)\ldots )$$
with $2r_0+q_0+2r_1+q_1+2r_2+q_2+\ldots +2r_{\ell}+q_{\ell}\leq p$
for some elements $\alpha_i\in \mathfrak{C}(F)$
for $i=0,\ldots ,\ell-1$
and polynomials $Q_i\in \mathsf{k}[z,w]$ of homogeneous degree
$q_i$ for $i=0,1,\ldots \ell$ and $\ell \geq 0$.
 
This obviously defines a filtered subalgebra of $End_{\mathsf{k}}(A)$.
It is easy to see that $D_0$ consists of endomorphisms defined by
$$\beta \mapsto a \beta $$
for some $a\in \mathsf{k}\mathfrak{C}(F)$. This follows 
because $\chi_{Q_0}$ is defined by multiplication with a constant
in $\mathsf{k}$. Note that 
$D_0=\mathcal{D}_0$. 

\vskip .1in

The theorem of section 5 implies:

\vskip .1in

\noindent \textbf{Theorem.} \textit{The $\star $-bracket on
$\mathsf{k}\mathfrak{C}(F)$ defined by the Kauffman bracket is
$D$-differentiable with respect to the filtered subalgebra $D$ defined
above.} 
 
\vskip .1in

Note that in this case the restriction of $\lambda_k$ to the second
variable for fixed $\alpha $ is given by
$$\varphi_0\chi_{P_k}+\varphi_1\chi_{P_{k-2}}+\ldots $$
with the polynomials $P_j$ from the theorem in section 5.  

\vskip .1in

Finally we will show that the restriction of $\lambda_2$ 
to one of the variables is \textit{not} a usual 
differential operator of order
$\leq 2$. 
In order to see this recall that 
$$\lambda_2(\alpha ,\beta )=\varphi_0\chi_{P_2}(\alpha
\triangleright \beta)+\varphi_1(\alpha \triangleright \beta),$$
because $P_0=1$. 
Even though \textit{only} the sum of the two terms is a well-defined
pairing $\mathsf{k}\mathfrak{C}(F)\otimes \mathsf{k}\mathfrak{C}(F)
\rightarrow \mathsf{k}\mathfrak{C}(F)$, it can still be checked whether
the differentiability formula holds separately for each term. 
But be aware that the value of each term depends on the choice of
diagram $\alpha \triangleright \beta $.

We want to check whether
$$
\beta \mapsto \lambda_2(\alpha \beta,\gamma)-\alpha \lambda_2(\beta,\gamma)
$$
is an operator of order $\leq 1$.
Let $P:=P_2$. The term $\varphi_0\chi_P$ is a differential operator of order $\leq2$.
Consider
$$
\Delta := \beta\mapsto \varphi_0 \chi_P[\alpha \triangleright \beta \triangleright\gamma ,
(\alpha \cap \gamma )\cup (\beta \cap \gamma)]
-\varphi_0\left(\alpha \cdot \chi_P[\beta \triangleright \gamma,
  \beta \cap \gamma]\right).
$$ 
Note that $\chi_P$ is defined by state summations over pairs of
crossings, so the above difference is determined by those states with at
least one of the two crossings on $\alpha $.  The application of
$\varphi_0$ does not change the formula since both $\alpha $ and
each term in $\chi_P[\beta \triangleright \gamma, \beta \cap
\gamma]$ have no crossings, and
$$\varphi_0(\alpha \triangleright \beta )=\varphi_0(\alpha )\varphi_0(\beta )$$
for all $\alpha, \beta \in \mathfrak{C}(F)$. Note that such a 
multiplicative property does \textit{not} hold for the higher order
projections $\varphi_i$ with $i>0$.  It follows that for all $\beta'$
$$ 
 \Delta (\beta \beta') - \beta \Delta \beta'
$$
involves only smoothings of pairs of crossings with one crossing in
$\alpha$ and one in $\beta$ and therefore is a multiple of $\beta'$.
Thus the first term of $\lambda_2$ has the differentiability
property of a differential operator of order $\leq 2$.

We now study the terms derived from the second contribution 
$\beta \mapsto \varphi_1(\alpha \triangleright \beta)$
thus whether
$$\delta \co \beta \mapsto \varphi_1(\alpha \triangleright \beta
\triangleright \gamma )-\alpha
\varphi_1(\beta \triangleright \gamma )$$
is an operator of order $\leq 1$ for all $\alpha, \gamma $.
This is the case if
$$\beta'\mapsto \delta (\beta \beta ')-\beta \delta (\beta ')$$
is an operator of order $0$ thus given by multiplying $\beta '$ by 
some element of $\mathsf{k}\mathfrak{C}(F)$.
If we let $\gamma =\emptyset $ and write the condition explicitly we get
$$\beta '\mapsto \varphi_1(\alpha \triangleright \beta
\triangleright \beta ')-\alpha \varphi_1(\beta \triangleright
\beta')-\beta \varphi_1(\alpha \triangleright \beta')-\beta \alpha
\varphi_1(\beta').$$
Thus $\beta '=\emptyset $ maps to some element
$\phi_1(\alpha \triangleright \beta )\in \mathsf{k}\mathfrak{C}(F)$.

Now consider the situation of $\alpha ,\beta$ two curves with no
crossings but both $\alpha $ and $\beta $ 
such that each state smoothing on $\beta \triangleright \beta '$
and $\alpha \triangleright \beta '$ does not involve an inessential
component while there exists a smoothing of $\alpha \triangleright 
\beta \triangleright \beta'$ involving an inessential component. 

In this case the differentiability condition is equivalent to 
$$\varphi_1(\alpha \triangleright \beta \triangleright \beta')=\alpha
\varphi_1(\beta \triangleright \beta ')+\beta \varphi_1(\alpha
\triangleright \beta ').$$
Then $\varphi_1(\beta ')=\varphi_1(\alpha \triangleright \beta
)=\varphi_1(\beta \triangleright \beta ')=\varphi_1(\alpha \triangleright
\beta ')=0$ while $\varphi_1(\alpha \triangleright \beta \triangleright
\beta ')\neq 0$.
Therefore the second term does not satisfy the differentiability
condition.  It follows that $\lambda_2$ restricted to one of the
variables is not a differential operator of order $\leq 2$.

\vskip .2in

\noindent \textbf{Remark.} The arguments above generalize to show that
the top term $\varphi_0\chi_{P_k}$ of $\lambda_k$ satisfies the 
condition of being a differentiable operator of order $\leq k$ in 
each variable. Note that this assertion is not precise in this form 
since we are
discussing homomorphisms from the module of diagrams into the algebra
$\mathsf{k}\mathfrak{C}(F)$. In fact, the operation $\triangleright$
used above is not well-defined but just a notation for a collection
of all diagrams. What we mean that the differentiability 
formula holds if we calculate 
$\varphi_k\chi_{P_k}$ on any diagram $\alpha \triangleright \beta$,
and multiplication in $\mathsf{k}\mathfrak{C}(F)$ is \textit{lifted}
to $\mathsf{k}\mathfrak{D}(F)$ in this way.

\vskip .1in

It seems a very difficult problem to compare the combinatorial filtered
subalgebra $D$ with the filtered subalgebra of differential operators
$\mathcal{D}$. But this problem is at the heart of relating the
combinatorial deformation with the geometric deformations of the
character variety.

\end{document}